\documentclass[12pt,twoside]{article}
\usepackage{graphicx}
\usepackage{amsfonts}
\usepackage{amsbsy}
\usepackage{amssymb}
\usepackage{amsmath}
\usepackage{cite}
\usepackage{pst-all}
\usepackage{pstricks}
\usepackage{graphics}
\def\bpsp{\begin{pspicture}}
\def\epsp{\end{pspicture}}

\newtheorem{theorem}{Theorem}[section]
\newtheorem{remark}[theorem]{Remark}
\newtheorem{example}[theorem]{Example}
\newtheorem{lemma}[theorem]{Lemma}
\newtheorem{corollary}[theorem]{Corollary}
\newtheorem{definition}[theorem]{Definition}
\newtheorem{proposition}[theorem]{Proposition}
\newtheorem{problem}[theorem]{Problem}
\newtheorem{note}{Note}
\newtheorem{case}{Case}
\newtheorem{conjecture}{Conjecture}
\newtheorem{question}{Question}

\newcommand{\bea}{\begin{eqnarray}}
\newcommand{\eea}{\end{eqnarray}}
\newcommand{\beq}{\begin{eqnarray*}}
\newcommand{\eeq}{\end{eqnarray*}}

\def\m4{\mbox{\rm ~(mod $4$)}}

\def \bd{\begin{definition}}
\def \ed{\end{definition}}
\def \bqu{\begin{question}}
\def \equ{\end{question}}
\def \bcc{\begin{conjecture}}
\def \ecc{\end{conjecture}}
\def \bt{\begin{theorem}}
\def \et{\end{theorem}}
\def \bl{\begin{lemma}}
\def \el{\end{lemma}}
\def \bc{\begin{corollary}}
\def \ec{\end{corollary}}
\def \be{\begin{equation}}
\def \ee{\end{equation}}
\def \ben{\begin{enumerate}}
\def \een{\end{enumerate}}
\def \ba{\begin{array}}
\def \ea{\end{array}}
\def \bp{\begin{proposition}}
\def \ep{\end{proposition}}
\def \bx{\begin{example}}
\def \ex{\end{example}}
\def \ex{\end{problem}}
\def \br{\begin{remark}}
\def \er{\end{remark}}
\def \bdsc{\begin{description}}
\def \edsc{\end{description}}

\def \bn{\begin{case}}
\def \en{\end{case}}
\def \bnt{\begin{note}}
\def \ent{\end{note}}
\def\1{1\!\!1}

\def\mm2{\mbox{\rm ~(mod $2$)}}
\def\m4{\mbox{\rm ~(mod $4$)}}

\def\qed{\nolinebreak\hfill\rule{.2cm}{.2cm}\par\addvspace{.5cm}}

\def\m{\mu}

\def\1{\textbf{1}}
\def\0{\textbf{0}}

\begin{document}
\title{On spectral spread of generalized distance matrix of a graph}
\author{ Hilal A. Ganie$^a$, S. Pirzada$^b$, A. Alhevaz$^c$, M. Baghipur$^d$\\
$^{a,b}${\em  Department of Mathematics, University of Kashmir, Srinagar, India}\\
$^{c,d}${\em Faculty of Mathematical Sciences, Shahrood University of Technology,}\\{\em P.O. Box: 316-3619995161, Shahrood, Iran}\\
$^a$hilahmad1119kt@gmail.com,~~$^b$pirzadasd@kashmiruniversity.ac.in\\
$^c$a.alhevaz@gmail.com, ~~ $^d$maryamb8989@gmail.com}
\date{}

\date{}
\pagestyle{myheadings} \markboth{Hilal, Pirzada, Alhevaz, Maryam}{On spectral spread of generalized distance matrix of a graph} \maketitle
\vskip 5mm
\noindent{\footnotesize \bf Abstract.} For a simple connected graph $G$, let $D(G)$, $Tr(G)$, $D^{L}(G)$ and $D^{Q}(G)$, respectively be the distance matrix, the diagonal matrix of the vertex transmissions, distance Laplacian matrix and the distance signless Laplacian matrix of a graph $G$. The convex linear combinations $D_{\alpha}(G)$ of $Tr(G)$ and $D(G)$ is defined as $D_{\alpha}(G)=\alpha Tr(G)+(1-\alpha)D(G)$,  $0\leq \alpha\leq 1$. As $D_{0}(G)=D(G), ~~~ 2D_{\frac{1}{2}}(G)=D^{Q}(G), ~~~ D_{1}(G)=Tr(G)$ and $D_{\alpha}(G)-D_{\beta}(G)=(\alpha-\beta)D^{L}(G)$, this matrix reduces to merging the distance spectral, distance Laplacian spectral and distance signless Laplacian spectral theories. Let $\partial_{1}(G)\geq \partial_{2}(G)\geq \dots \geq \partial_{n}(G)$ be the eigenvalues of $D_{\alpha}(G)$ and let $D_{\alpha}S(G)=\partial_{1}(G)-\partial_{n}(G)$ be the generalized distance spectral spread of the graph $G$. In this paper, we obtain some bounds for the generalized distance spectral spread $D_{\alpha}(G)$. We also obtain relation between the generalized distance spectral spread $D_{\alpha}(G)$ and the distance spectral spread $S_{D}(G)$.  Further, we obtain the lower bounds for $D_{\alpha}S(G)$ of bipartite graphs involving different graph parameters and we characterize the extremal graphs for some cases. We also obtain lower bounds for $D_{\alpha}S(G)$ in terms of clique number and independence number of the graph $G$  and characterize the extremal graphs for some cases.
\vskip 3mm

\noindent{\footnotesize Keywords: Distance spectrum (matrix); Distance signless Laplacian spectrum (matrix); transmission regular graph; generalized distance matrix, generalized distance spectral spread.}

\vskip 3mm
\noindent {\footnotesize AMS subject classification: 05C50, 05C12, 15A18.}

\section{Introduction}

All graphs considered in this paper are finite, undirected and simple. Let $G$ be such a graph with vertex set $V(G)$ and edge set $E(G)$, denoted by $G=(V(G),E(G))$. We use standard terminology; for concepts not defined here, we refer the reader to any standard graph theory monograph, such as \cite{BH} or \cite{SP}.\\
\indent In $G$, the \textit{distance} between two vertices $u,v\in V(G),$ denoted by $d_{uv}$, is defined as the length of a shortest path between $u$ and $v$. The \textit{diameter} of $G$ is the maximum distance between any two vertices of $G$. The \textit{distance matrix} of $G$, denoted by $D(G)$, is defined as $D(G)=(d_{uv})_{u,v\in V(G)}$. Till now, the distance spectrum of a connected graph has been investigated extensively. The \textit{transmission} $Tr_{G}(v)$ of a vertex $v$ is defined to be the sum of the distances from $v$ to all other vertices in $G$, i.e., $Tr_{G}(v)=\sum\limits_{u\in V(G)}d_{uv}$. A graph $G$ is said to be $k$-\textit{transmission regular} if $ Tr_{G}(v)=k,$ for each $ v\in V(G). $  The \textit{Wiener index} of a graph $ G, $ denoted by $W(G), $ is the sum of distances between all unordered pairs of vertices in $G$. Clearly, $W(G)=\frac{1}{2}\displaystyle\sum_{v\in V(G)}Tr_{G}(v)$. For a graph $G$ with $ V(G)=\{v_{1},v_{2},\ldots,v_{n}\}, $  $Tr_G(v_i)$ has been referred as the \textit{transmission degree} $Tr_{i}$  and hence the \textit{transmission degree sequence} is given by $\{Tr_{1},Tr_{2},\ldots,Tr_{n}\}$. The second transmission degree of $v_{i}$, denoted by $T_{i}$ is given by $T_{i}=\displaystyle\sum_{j=1}^{n}d_{ij}Tr_{j}$. Let $Tr(G)=diag(Tr_1,Tr_2,\ldots,Tr_n)$ be the diagonal matrix of vertex transmissions of $G$. Aouchiche and  Hansen \cite{AH1,AH2} introduced the Laplacian and the signless Laplacian for the distance matrix of a connected graph.
The matrix $D^L(G)=Tr(G)-D(G) $ is called the \textit{distance Laplacian matrix} of $G$, while the matrix  $D^{Q}(G)=Tr(G)+D(G)$ is called the \textit{distance signless Laplacian matrix of $G$.} For some recent research we refer to \cite{amh,abhp,ap,xiom} and the references therein.

Recently in  \cite{CHT}, Cui et al. introduced the \textit{generalized distance matrix} $D_{\alpha}(G)$ as a convex combinations of $Tr(G)$ and $D(G)$, defined as $D_{\alpha}(G)=\alpha Tr(G)+(1-\alpha)D(G)$,  for $0\leq \alpha\leq 1$. Since $D_{0}(G)=D(G), ~~~ 2D_{\frac{1}{2}}(G)=D^{Q}(G), ~~~ D_{1}(G)=Tr(G)$ and  $D_{\alpha}(G)-D_{\beta}(G)=(\alpha-\beta)D^{L}(G)$, any result regarding the spectral properties of generalized distance matrix, has its counterpart for each of these particular graph matrices, and these counterparts follow immediately from a single proof. In fact, this matrix reduces to merging the distance spectral, distance Laplacian spectral and distance signless Laplacian spectral theories.

Since the matrix $D_{\alpha}(G)$ is a real symmetric matrix, therefore its eigenvalues can be
arranged as $\partial_{1}\geq \partial_{2}\geq \cdots \geq \partial_{n}$, where the largest eigenvalue $\partial_{1}$ is called the \textit{generalized distance spectral radius} of $G$. (From now onwards, we will denote $\partial_1(G)$ by $\partial_1$). As $D_{\alpha}(G)$ is non-negative and irreducible, by the Perron-Frobenius theorem, $ \partial_1$ is the unique eigenvalue and there is a unique positive unit eigenvector $X$ corresponding to $ \partial_1,$ which is called the \textit{generalized distance Perron vector} of $G$.

Let $M$ be a symmetric matrix of order $n$ having eigenvalues $w_i, i=1,2,\dots,n$. If $w_1$ and $w_n$ are respectively the largest and smallest eigenvalues of the symmetric matrix $M$, then the spread of the matrix $M$ is defined as $s(M)=w_1-w_n$. In the literature, several papers can be seen regarding the parameter $s(M)$ for any symmetric matrix $M$. The spread of a matrix is a very attractive topic and as such the investigation of the spread of some matrices of a graph become interesting. When $M$ is restricted to a particular graph matrix, the parameter $s(M)$ has attracted much attention of the researchers as is clear from the fact that  various papers can be found in the literature in this direction. For a particular graph matrix (like adjacency, Laplacian, signless Laplacian etc), the much studied problem about the parameter $s(M)$ is to obtain bounds in terms of various graph parameters. Another problem worth to mention is to characterize the extremal graphs for the parameter $s(M)$ for a graph matrix, in some class of graphs.\\
\indent Let $\rho_1\geq \rho_2\geq \dots\rho_n$ be the distance eigenvalues of the graph $G$. The distance spread of a connected graph $G$ was considered in \cite{yzlws} and is defined as $S_{D}(G)=\rho_1-\rho_n$. Various lower and upper bounds for the parameter $S_{D}(G)$ can be found in in \cite{yzlws}.\\
\indent Motivated by the results obtained for the distance spread, You et al. \cite{you} put forward concept of the distance signless Laplacian spread of a connected graph $G$ as $ S_{D^{Q}}(G)=q^{D}_1-q^{D}_n$ and obtained various lower and upper bounds for the parameter $S_{D}(G)$ in terms of different graph parameters.

We define the generalized distance spectral spread of a graph $G$ as the difference between the
largest and smallest eigenvalues of $D_{\alpha}(G)$, that is,
\begin{align*}
 D_{\alpha}S(G)=\partial_{1}(G)-\partial_{n}(G).
\end{align*}
It is clear from the above discussion that $ D_{0}S(G)=S_{D}(G)$ and $2D_{\frac{1}{2}}S(G)=S_{D^{Q}}(G)$. Therefore, the parameter $D_{\alpha}S(G)$ is a generalization of the already studied parameters $S_{D}(G)$ and $S_{D^{Q}}(G)$. The motive of this paper is to obtain bounds for the parameter  $D_{\alpha}S(G)$, in terms of different graph parameters and to characterize the extremal graphs. \\
\indent The rest of the paper is organized as follows.  In Section 2, we obtain some bounds for the generalized distance spectral spread $D_{\alpha}(G)$ of graphs. We also obtain relation between the generalized distance spectral spread $D_{\alpha}(G)$ and the distance spectral spread $S_{D}(G)$.  In Section 3, we obtain lower bounds for the generalized distance spectral spread $D_{\alpha}(G)$ of bipartite graphs involving different graph parameters and characterize the extremal graphs for some cases. In Section 4, we obtain lower bounds for the generalized distance spectral spread $D_{\alpha}(G)$ in terms of clique number and independence number of the graph $G$  and characterize the extremal graphs for some cases.

\section{Bounds for $D_{\alpha}(G)$ for any graph}

In this section, we obtain some bounds for $D_{\alpha}S(G)$, in terms of various graph parameters. We also establish a relation between the generalized distance spectral spread and the distance spread of a graph $G$. We start with known results which will be used in the proofs of our main results in the sequel.

\begin{lemma}{\em  \cite{CHT}} \label{lem 2.1} Let $G$ be a connected graph of order $n$. Then,
$\partial_{1} (G)\geq \frac{2W(G)}{n},$ with equality if and only if $G$ is a transmission regular graph.
\end{lemma}

\indent The following Lemma can be found in \cite{w}.\\

 \begin{lemma}
  Let $X$ and $Y$ be Hermitian matrices of order $n$ such that $Z=X+Y$. Then
  \begin{align*}
  \lambda_k(Z)\leq \lambda_j(X)+\lambda_{k-j+1}(Y),~~ n\geq k\geq j\geq 1,\\
  \lambda_k(Z)\geq \lambda_j(X)+\lambda_{k-j+n}(Y), ~~n\geq j\geq k\geq 1,
  \end{align*}
 where $\lambda_i(M)$ is the $i^{th}$ largest eigenvalue of the matrix $M$. In either of these inequalities, equality holds if and only if there exists a unit vector that is an eigenvector to each of the three eigenvalues involved.\\
\end{lemma}

The proof of the following lemma is similar to that of Lemma 2 in \cite{I}, so is omitted here.

\begin{lemma}\label{lemm2.3} A connected graph $G$ has two distinct $D_{\alpha}(G)$ eigenvalues if and only if $G$ is a complete graph.
\end{lemma}

If $G$ is a $k$-transmission regular graph, then $Tr(G)=k I_n$ (where $I_n$ is the identity matrix of order $n$) and so
\begin{align*}
D_{\alpha}(G)=\alpha Tr(G)+(1-\alpha)D(G)= k\alpha I_n +(1-\alpha)D(G).
\end{align*}
Therefore, if $\rho_1\geq \rho_2\cdots\geq \rho_n$ are the distance eigenvalues of the $k$-transmission regular graph $G$, its generalized distance eigenvalues are $\partial_{i}=k\alpha+(1-\alpha)\rho_i$, $1\leq i\leq n$. Thus, for a  $k$-transmission regular graph $G$, we have
\begin{align*}
D_{\alpha}S(G)=\partial_{1}-\partial_{n}=(1-\alpha)(\rho_1-\rho_n)=(1-\alpha)S_D(G).
\end{align*}
Hence, for a $k$-transmission regular graph $G$, the study of generalized distance spectral spread $D_{\alpha}S(G)$ is same as distance spread $S_D(G)$.\\
\indent For $\frac{1}{2}\leq \alpha\leq 1$, it has been shown \cite{CHT} that the generalized distance matrix $D_{\alpha}(G)$ of the connected graph $G$ is a positive semi-definite matrix. Therefore, for $\frac{1}{2}\leq \alpha\leq 1$, clearly $\partial_{n}\geq 0$ implies that $D_{\alpha}S(G)=\partial_{1}-\partial_{n}\leq \partial_{1}$, with equality if and only if $\partial_{n}=0$. From this, for $\frac{1}{2}\leq \alpha\leq 1$, it is clear that any upper bound for generalized distance spectral radius gives an upper bound for the generalized distance spectral spread $D_{\alpha}S(G)$.\\

\indent The following result gives a relation between the generalized distance spectral spread $D_{\alpha}S(G)$ and the distance spread  $S_{D}(G)$ of a connected graph $G$.

\begin{theorem}
Let $G$ be a connected graph of order $n$ having transmission degree sequence $\{Tr_{1},Tr_{2},\ldots,Tr_{n}\}$.  If $Tr_{\max}=\max_{1\leq i\leq n}Tr_i$ and $Tr_{\min}=\min_{1\leq i\leq n}Tr_i$, then
\begin{align}\label{r}
\Big|\alpha(Tr_{\max}-Tr_{\min})-(1-\alpha)S_{D}(G)\Big|\leq D_{\alpha}S(G)\leq \alpha\Big(Tr_{\max}-Tr_{\min}\Big)+(1-\alpha)S_{D}(G).
\end{align} Equality occurs on both sides if and only if $G$ is a transmission regular graph..
\end{theorem}
{\bf Proof.}  Let $G$ be a connected graph of order $n$ having generalized distance matrix $D_{\alpha}(G)$. Since
$D_{\alpha}(G)=\alpha Tr(G)+(1-\alpha)D(G)$, by taking $k=1$, $j=1$ in first inequality and $k=1$, $j=n$ in second inequality of Lemma 2.2, it follows that
\begin{align}
\lambda_n(\alpha Tr(G)) +\lambda_1((1-\alpha)D(G))\leq \lambda_1(D_{\alpha}(G))\leq \lambda_1(\alpha Tr(G)) +\lambda_1((1-\alpha)D(G)).
\end{align}
Again taking $k=n$, $j=1$ in first inequality and $k=n$, $j=n$ in second inequality of Lemma 2.2, we get
\begin{align}
\lambda_n(\alpha Tr(G)) +\lambda_n((1-\alpha)D(G))\leq \lambda_n(D_{\alpha}(G))\leq \lambda_1(\alpha Tr(G)) +\lambda_n((1-\alpha)D(G)).
\end{align}
From (2.2) and (2.3), we have
\begin{align*}
\lambda_1(D_{\alpha}(G))-\lambda_n(D_{\alpha}(G))\leq \alpha\Big(\lambda_1(Tr(G))-\lambda_n( Tr(G))\Big)+(1-\alpha)\Big(\lambda_1(D(G))-\lambda_n(D(G))\Big)
\end{align*}
and
\begin{align*}
\lambda_1(D_{\alpha}(G))-\lambda_n(D_{\alpha}(G))\geq \Big|\alpha\Big(\lambda_1(Tr(G))-\lambda_n( Tr(G))\Big)-(1-\alpha)\Big(\lambda_1(D(G))-\lambda_n(D(G))\Big)\Big|.
\end{align*}
The result now follows from these inequalities. Equality occurs in (2.1) if and only if equality occurs in (2.2) and (2.3). Suppose that equality occurs on the right of (2.2), then by Lemma 2.2, the eigenvalues $\lambda_1(D_{\alpha}(G)), \lambda_1( Tr(G))$ and $\lambda_1(D(G))$ of the matrices $D_{\alpha}(G), Tr(G)$ and $D(G)$  have the same unit eigenvector $\boldsymbol{x}$. Similarly, if equality occurs on the left of (2.2), then again by Lemma 2.2, the eigenvalues $\lambda_1(D_{\alpha}(G)), \lambda_n( Tr(G))$ and $\lambda_1(D(G))$ of the matrices $D_{\alpha}(G), Tr(G)$ and $D(G)$ have the same unit eigenvector $\boldsymbol{x}$. From this it follows that the  eigenvalues $\lambda_1( Tr(G))$ and $\lambda_n( Tr(G))$ of the matrix $Tr(G)$ have the same eigenvector, which is only possible if $\lambda_1( Tr(G))=\lambda_n( Tr(G))$. This gives that $Tr(G)=\lambda_1( Tr(G))I_n$, that is, $G$ is a $\lambda_1( Tr(G))$ transmission regular graph. In a similar way we can discuss the equality in (2.3).\\
\indent Conversely, if $G$ is a $k$-transmission regular graph then $D_{\alpha}(G)= k\alpha I_n +(1-\alpha)D(G)$ and so equality holds in (2.2) and (2.3).\qed

From inequality (2.2), we obtain the following bounds for the generalized distance spectral radius $\partial_1$
\begin{align}
\alpha Tr_{\min} +(1-\alpha)\rho_1\leq \partial_1 \leq \alpha Tr_{\max} +(1-\alpha)\rho_1,
\end{align}
with equality occurs on both sides if and only if $G$ is a transmission regular graph. In the literature, we can find various bounds for the distance spectral radius $\rho_1$. All those bounds together with inequality (2.4) give bounds for the generalized distance spectral radius $\partial_1$.\\
\indent  Similarly, from inequality (2.3), we get the following bounds for the smallest generalized distance eigenvalue $\partial_n$
\begin{align}
\alpha Tr_{\min} +(1-\alpha)\rho_n\leq \partial_n \leq \alpha Tr_{\max} +(1-\alpha)\rho_n,
\end{align}
with equality on both sides if and only if $G$ is a transmission regular graph.

\indent The following result gives the lower bound for $D_{\alpha}S(G)$ in terms of Weiner index $W$ and the generalized distance spectral radius $\partial_1$.

\begin{theorem} If $G$ is a connected  graph of order $n$, then
\begin{eqnarray}\label{e17}
D_{\alpha}S(G)\geq \frac{n}{n-1}\partial_{1}-\frac{2\alpha W(G)}{n-1},
\end{eqnarray}
equality holds if and only if $G\cong K_{n}.$
\end{theorem}
{\bf Proof.} Let $\partial_{1},\ldots,\partial_{n}$ be the generalized distance eigenvalues of $G.$ Note that $\sum_{i=1}^{n}\partial_{i}=2\alpha W(G).$ So
$2\alpha W(G)\geq \partial_{1}+(n-1)\partial_{n}$ and therefore $\partial_{n}\leq \frac{(2\alpha W-\partial_{1})}{n-1}.$ Thus
\begin{eqnarray*}
D_{\alpha}S(G)\geq \partial_{1}-\frac{(2\alpha W(G)-\partial_{1})}{n-1}=\frac{n}{n-1}\partial_{1}-\frac{2\alpha W(G)}{n-1}.
\end{eqnarray*}
Equality in (\ref{e17}) holds if and only if $\partial_{2}=\cdots=\partial_{n}$.  That is, if and only if $G$ has exactly two distinct eigenvalues. Using  Lemma \ref{lemm2.3}, it follows that $G$ is the complete graph $K_{n}.$\\
\indent Conversely, if $G\cong K_{n}$, it can be seen by direct computations that equality occurs in (\ref{e17}).  \qed

\indent The following is another lower bound for $D_{\alpha}S(G)$ in terms of Weiner index $W$ and the generalized distance spectral radius $\partial_1$.

\begin{theorem} Let $G$ be a connected graph of order $n$. Then
\begin{eqnarray*}
D_{\alpha}S(G)\geq \partial_{1}-\sqrt{\frac{2(1-\alpha)^2\sum_{1\leq i<j\leq n}(d_{ij})^{2}+\alpha ^2\sum_{i=1}^{n}Tr^{2}_{i}-\partial^{2}_{1}}{n-1}},
\end{eqnarray*}
equality holds if and only if $G\cong K_{n}$.
\end{theorem}
{\bf Proof.} Note that
\begin{eqnarray}\label{e18}
2(1-\alpha)^2\sum_{1\leq i<j\leq n}(d_{ij})^{2}+\alpha ^2\sum_{i=1}^{n}Tr^{2}_{i}=\sum_{i=1}^{n}\partial^{2}_{i}\geq \partial^{2}_{1}+(n-1)\partial^{2}_{n}.
\end{eqnarray}
Therefore, $$\partial_{n}\leq \sqrt{\frac{2(1-\alpha)^2\sum_{1\leq i<j\leq n}(d_{ij})^{2}+\alpha ^2\sum_{i=1}^{n}Tr^{2}_{i}-\partial^{2}_{1}}{n-1}}.$$ Thus $$D_{\alpha}S(G)\geq \partial_{1}-\sqrt{\frac{2(1-\alpha)^2\sum_{1\leq i<j\leq n}(d_{ij})^{2}+\alpha ^2\sum_{i=1}^{n}Tr^{2}_{i}-\partial^{2}_{1}}{n-1}}.$$
Equality in (\ref{e18}) holds if and only if $ \partial_{2}=\cdots=\partial_{n}.$  That is, if and only  $ G $ has exactly two distinct eigenvalues. Using Lemma \ref{lemm2.3}, it follows that $G$ is the complete graph $K_{n}.$ \\
\indent Conversely, if $G\cong K_{n}$, it can be seen by direct computations that equality occurs in (\ref{e18}).\qed

The following observation is immediate from Lemma 2.1 and Theorem 2.6.

\begin{corollary} Let $G$ be a connected graph of order $n$ having Wiener index $W(G)$. Then
\begin{eqnarray*}
D_{\alpha}S(G)\geq\frac{1}{n}\left(2W(G)-\sqrt{\frac{n^{2}(2(1-\alpha)^2\sum_{1\leq i<j\leq n}(d_{ij})^{2}+\alpha ^2\sum_{i=1}^{n}Tr^{2}_{i})-4W^{2}(G)}{n-1}}\right),
\end{eqnarray*}
with equality if and only if $G\cong K_{n}$.
\end{corollary}

\indent Now, we obtain a lower bound for the generalized distance spectral spread $D_{\alpha}S(G)$ in terms of Weiner index $W$.
\begin{theorem} Let $G$ be a connected graph of order $n$ having Wiener index $W(G)$. Then
\begin{eqnarray*}
D_{\alpha}S(G)\geq \frac{2}{n}\sqrt{n\left(2(1-\alpha)^2\sum_{1\leq i< j\leq n}(d_{ij})^{2}+\alpha ^2\sum_{i=1}^{n}Tr_{i}^{2}\right)-4\alpha ^2W^{2}(G)}.
\end{eqnarray*}
\end{theorem}
{\bf Proof.} Let $\partial_{1}\geq \partial_{2}\geq \ldots\geq \partial_{n}$ be the generalized distance eigenvalues of $G.$ We assume that the non-negative real numbers $a_{2},a_{3},\ldots,a_{n-1}, $ $ b_{2},b_{3},\ldots,b_{n-1}$ are defined by the equations $\partial_{i}^{2}=a_{i}\partial_{1}^{2}+b_{i}\partial^{2}_{n},\quad a_{i}+b_{i}=1,\quad i=2,3,\ldots ,n-1$. Putting $a=1+\sum_{i=2}^{n-1}a_{i}$ and $b=1+\sum_{i=2}^{n-1}b_{i},$ then $a+b=n$ and $\sum_{i=1}^{n}\partial_{i}^{2}=a\partial_{1}^{2}+b\partial_{n}^{2}$. On the other hand, we have $\partial_{i}^{2}=(a_{i}\partial_{1}^{2}+b_{i}\partial_{n}^{2})(a_{i}+b_{i}),$ so that $\partial_{i}\geq a_{i}\partial_{1}+b_{i}\partial_{n}, \quad i=2,3,\ldots, n-1$, and therefore $\sum_{i=1}^{n}\partial_{i}\geq a\partial_{1}+b\partial_{n}$. Now, observing that $ab\leq \frac{(a+b)^{2}}{4}=\frac{n^{2}}{4},$ we obtain
\begin{eqnarray}\label{e19}
n\sum_{i=1}^{n}\partial_{i}^{2}-\left(\sum_{i=1}^{n}\partial_{i}\right)^{2}&\leq & (a\partial_{1}^{2}+b\partial_{n}^{2})(a+b)-(a\partial_{1}+b\partial_{n})^{2}\nonumber\\
&=&ab(\partial_{1}-\partial_{n})^{2}\leq \frac{n^{2}}{4}(\partial_{1}-\partial_{n})^{2}.
\end{eqnarray}
Consequently, from (\ref{e19}),  we get
\begin{eqnarray*}
D_{\alpha}S(G)\geq \frac{2}{n}\sqrt{n\left(2(1-\alpha)^2\sum_{1\leq i< j\leq n}(d_{ij})^{2}+\alpha ^2\sum_{i=1}^{n}Tr_{i}^{2}\right)-4\alpha ^2W^{2}(G)},
\end{eqnarray*}
as desired.\qed

For a real rectangular matrix $M=(w_{ij})_{mxn}$, let $\|M\|_F=\left(\sum\limits_{i=1}^{m}\sum\limits_{j=1}^{n}|w_{ij}|^2\right)^\frac{1}{2}$ be the Frobenius norm of $M$. If $M$ is a square matrix, its trace is denoted by
$tr M$. If $w_1$ and $w_n$ are respectively the largest and smallest eigenvalues of the symmetric matrix $M$, the spread of the matrix $M$ is defined as $s(M)=w_1-w_n$. The following observation is due to Mirsky \cite{mirk}.
\begin{lemma}
Let $M$ be an $n$-square normal matrix. Then
\begin{align*}
s(M)\leq \Big(2\|M\|_F^{2}-\frac{2}{n}(tr M)^2\Big)^{\frac{1}{2}},
\end{align*}
\end{lemma}
equality occurs if and only if the eigenvalues $w_1,w_2,\dots,w_n$ of $M$ satisfy the
$w_2=w_3=\dots=w_{n-1}=\frac{w_1+w_n}{2}$.

\begin{theorem}
Let $G$ be a connected graph of order $n$ having Wiener index $W$. Then
\begin{align*}
D_{\alpha}S(G)\leq \Big(2(1-\alpha)^2\sum\limits_{i\ne j}(d_{ij})^2+\alpha^2 \sum\limits_{i=1}^{n}Tr_i^{2}-\frac{8}{n}\alpha^2W^2\Big)^{\frac{1}{2}},
\end{align*} with equality if and only if $\partial_2=\partial_3=\cdots=\partial_{n-1}=\frac{\partial_1+\partial_n}{2}$.
\end{theorem}
{\bf Proof.} Since the matrix $D_{\alpha}(G)$ is a normal matrix with $\|D_{\alpha}(G)\|_F^{2}=(1-\alpha)^2\sum\limits_{i\ne j}(d_{ij})^2+\alpha^2 \sum\limits_{i=1}^{n}Tr_i^{2}$ and $tr D_{\alpha}(G)=2\alpha W$, the result follows from Lemma 2.9, by taking $M=D_{\alpha}(G)$. \qed
\indent For $\alpha=1$, it can be seen that equality occurs in Theorem 2.10 if and only if $G\cong K_n$. In general, equality occurs in Theorem 2.10 if and only if $\partial_2=\partial_3=\dots=\partial_{n-1}=\frac{\partial_1+\partial_n}{2}$. This gives $tr D_{\alpha}(G)=n \partial_2$, which implies that $2\alpha W=n \partial_2$, that is, $\partial_2=\frac{2\alpha W}{n}$. Since
\begin{align}
\frac{2\alpha W}{n}=\alpha\frac{e^t D_{\alpha}(G) e }{e^t e}\leq \alpha \partial_1,
\end{align}
where $e$ is the all one $n$-vector, it follows that $\partial_2 \leq \alpha \partial_1$. So, we have $\partial_1+ \partial_n =2\partial_2\leq 2 \alpha \partial_1$, implying that $\partial_n\leq (2 \alpha -1)\partial_1 $, if $\alpha>\frac{1}{2}$ and $ (1-2 \alpha )\partial_1+\partial_n\leq 0$, if $\alpha\leq \frac{1}{2}$. Since equality occurs in (2.9) if and only if $G$ is a transmission regular graph, it follows that equality occurs in Theorem 2.10 if and only if $G$ is a transmission regular graph with $\partial_2=\partial_3=\dots=\partial_{n-1}=\alpha \partial_1$ and $\partial_n\leq (2 \alpha -1)\partial_1$, for $\alpha>\frac{1}{2}$ and $(1-2 \alpha )\partial_1+\partial_n\leq 0$, for $\alpha\leq \frac{1}{2}$. In particular, if $\alpha=0$, then equality occurs if and only if $\partial_1=- \partial_n,  \partial_2=\cdots=\partial_{n-1}=0$.

\section{\bf Bounds for  $D_{\alpha}S(G)$ for bipartite graphs}

In this section, we obtain lower bounds for $D_{\alpha}S(G)$ in terms of maximum degree $\Delta$, order $n$, transmission degrees, average distance degrees and  Weiner index of a bipartite graph $G$.

Let $M$ be a real matrix of order $n$ described in the following block form
\begin{align}
M =
\begin{pmatrix}\label{x}
M_{11} & M_{12}&\cdots & M_{1t}\\
\vdots & \vdots & \cdots & \vdots\\
M_{t1} & M_{t2}&\cdots & M_{tt}
\end{pmatrix},
\end{align} where the diagonal blocks $M_{ii}$ are $n_i\times n_i$ matrices for any $i \in \{1,2,\dots,t\}$ and
$n=n_1+\dots+n_t$. For any $i,j \in  \{1,2,\dots,t\}$, let $b_{ij}$ denote the average row sum of $M_{ij}$, that is, $b_{ij}$
is the sum of all entries in $M_{ij}$ divided by the number of rows. Then $B(M)=(b_{ij})$ is called the quotient matrix of $M$. If in addition for each pair $i,j$, $M_{ij}$ has constant
row sum, then $B(M)$ is called the equitable quotient matrix of $M$. The following  Lemma can be found in \cite{BH}.
\begin{lemma}
Let $M$ be a symmetric matrix which has the block form as (\ref{x}) and $B$ be
the quotient matrix of $M$. Then the eigenvalues of $B$ interlace the eigenvalues of $M$. That is, if $a_1\geq a_2\geq \cdots\geq a_n$ and $b_1\geq b_2\cdots\geq b_r$ are the eigenvalues of the matrices $M$ and $B$, respectively. Then $a_i\geq b_i\geq a_{n-r+i}$, for all $i=1,2,\dots,r$.
\end{lemma}
\indent The following observation can be found in \cite{cds}.
\begin{lemma}
Let $M$ be a matrix of order $n$ having eigenvalues $a_1\geq a_2\geq \cdots\geq a_n$ and let $B$ be a principle matrix of $M$ of order $r$ having eigenvalues $b_1\geq b_2\cdots\geq b_r$, then $a_i\geq b_i\geq a_{n-r+i}$, for all $i=1,2,\dots,r$.
\end{lemma}

\indent The following observation can be found in \cite{CHT}.
\begin{lemma}
Let $G$ be a connected graph of order $n$ and $S$ be a subset of the vertex set $V(G)$, such that $N(x)=N(y)$ for any $x, y\in S$, where $|S|=t$.\\
(i) If $S$ is an independent set, then $Tr(x)$ is a constant for each $x\in S$ and $D_{\alpha}(G)$ has $\alpha(Tr(x)+2)-2$ as an eigenvalue with multiplicity at least $t-1$.
(ii) If $S$ is a clique, then $Tr(x)$ is a constant for each $x\in S$ and $D_{\alpha}(G)$ has $\alpha(Tr(x)+1)-1$ as an eigenvalue with multiplicity at least $t-1$.
\end{lemma}

\indent The following result gives the generalized distance spectrum of the complete bipartite graph $K_{r,s}$.
\begin{lemma}
The generalized distance spectrum of complete bipartite graph $K_{r,s}$ consists of eigenvalue $\alpha(2r+s)-2$ with multiplicity $r-1$, eigenvalue $\alpha(2s+r)-2$ with multiplicity $s-1$ and the remaining two  eigenvalues as $x_1, x_2$, where $x_1,x_2=\frac{\alpha(s+r)+2(s+r)-4\pm \sqrt{(r^2+s^2)(\alpha-2)^2+2rs(\alpha^2-2)}}{2}$.
\end{lemma}
{\bf Proof.} Let $V_1$ and $V_2$ be the partite sets of $K_{r,s}$ such that the degree of each vertex in $V_1$ is $r$ and the degree of each vertex in $V_2$ is $s$. It is clear that $Tr(v_i)=2r+s-2$, for all $v_i\in V_1$ and $Tr(u_j)=2s+r-2$, for all $u_j\in V_2$. Now, using Lemma 3.3 with $S=V_1$ and then with $S=V_2$, we get the eigenvalue $\alpha(2r+s)-2$ with multiplicity $r-1$ and the eigenvalue $\alpha(2s+r)-2$ with multiplicity $s-1$. The remaining two eigenvalues (as the quotient matrix for the graph $K_{r,s}$ is an equitable quotient matrix) can be obtained from the quotient matrix of the matrix $D_{\alpha}(K_{r,s})$ given by $B =
\begin{pmatrix}\label{y}
\alpha s+2r-2 & s(1-\alpha) \\
r(1-\alpha) & \alpha r+2s-2
\end{pmatrix}$.  \qed

 \indent Let $t_{v_i}=\frac{1}{d_{G}(v_i)}\sum\limits_{v_j,v_iv_j\in E(G)}Tr(v_j)$ be the average distance degree of the vertex $v_i$. The following result gives a lower bound  for $D_{\alpha}S(G)$ in terms of maximum degree $\Delta$, order $n$, average distance degrees and  Weiner index $W$ of a bipartite graph $G$.
\begin{theorem}
Let $G$ be a connected bipartite graph on $n\geq 3$ vertices with maximum
degree $\Delta$ and Wiener index $W$. Suppose that $d_{G}(v_1)=d_{G}(v_2)=\cdots=d_{G}(v_k)=\Delta$. \\
{\bf (i)}. If $\Delta\leq n-2$, then
\begin{align*}
D_{\alpha}S(G)\geq \max_{1\leq i\leq k}\frac{\sqrt{\alpha_i^{2}-4\beta_i(\Delta+1)(n-\Delta-1)}}{(\Delta+1)(n-\Delta-1)},
\end{align*}  where $\alpha_i=\alpha n(t_{v_i}\Delta+Tr(v_i)-2\Delta^2)+2n\Delta^2+(\Delta+1)(2W-2t_{v_i}\Delta-2Tr(v_i))$ and $\beta_i=2\alpha W(t_{v_i}\Delta+Tr(v_i)-2\Delta^2)+4W\Delta^2-(t_{v_i}\Delta+Tr(v_i))^2$.\\
{\bf (ii)}. If $\Delta= n-1$, then $D_{\alpha}S(G)=\left\{\begin{array}{lr}n+\sqrt{n^2-3n+3}, &\mbox{if $\alpha =0$,}\\
\sqrt{(\alpha-2)^2(n^2-2n+2)+2(n-1)(\alpha^2-2)}, & \mbox{ if $\alpha\ne 0$.}
\end{array}\right.$
\end{theorem}
{\bf Proof.} Let $G$ be a connected bipartite graph on $n\geq 3$ vertices with maximum
degree $\Delta$, Wiener index $W$ and  $d_{G}(v_1)=d_{G}(v_2)=\cdots=d_{G}(v_k)=\Delta$.\\
{\bf (i).} Let  $\Delta\leq n-2$. For $1\leq i\leq k$, let $N(v_i)=\{u_{i_1},u_{i_2},\dots,u_{i_\Delta}\}$ be the neighbor set and $N[v_i]=N(v_i)\cup\{v_i\}$. Let $V(G)=\{v_i,u_{i_1},u_{i_2},\dots,u_{i_\Delta},w_1,w_2,\dots,w_{n-\Delta-1}\}$ and let $V(G)=V_1\cup V_2$ be the bipartite partition of $V(G)$. Since $G$ is bipartite, if $v_i\in V_1$, then $u_{ij}\in V_2$ and if $v_i\in V_2$, then $u_{ij}\in V_1$, for $j=1,2,\dots,\Delta$. By suitably labelling the vertices of $G$, it can be seen that the generalized distance matrix of the graph $G$ can be written as
$D_{\alpha}(G) =
\begin{pmatrix}\label{y}
X & (1-\alpha) Y\\
(1-\alpha) Y^t & Z
\end{pmatrix},$ where
\begin{align*}
X =\begin{pmatrix}
\alpha Tr(v_i) & 1-\alpha& 1-\alpha&\cdots & 1-\alpha\\
1-\alpha & \alpha Tr(u_{i_1}) & 2(1-\alpha) & \cdots & 2(1-\alpha)\\
\vdots & \vdots & \vdots & \cdots & \vdots\\
1-\alpha &2(1-\alpha)  & 2(1-\alpha) & \cdots &\alpha Tr(u_{i_\Delta})
\end{pmatrix},
Y =\begin{pmatrix}
d(v_i,w_1)  &\cdots & d(v_i,w_{n-\Delta-1})\\
d(u_{i_1},w_1) & \cdots & d(u_{i_1}w_{n-\Delta-1})\\
\vdots  & \cdots & \vdots\\
 d(u_{i_{\Delta}},w_1)&\cdots & d(u_{i_{\Delta}}w_{n-\Delta-1})
\end{pmatrix}
\end{align*}
and $Z=\begin{pmatrix}
\alpha Tr(w_1)& (1-\alpha)d(w_1,w_2)  &\cdots & (1-\alpha)d(w_1,w_{n-\Delta-1})\\
(1-\alpha)d(w_2,w_1)& \alpha Tr(w_2)  & \cdots & (1-\alpha)d(w_2,w_{n-\Delta-1})\\
\vdots  & \vdots& \cdots & \vdots\\
 (1-\alpha)d(w_{n-\Delta-1},w_1)& (1-\alpha)d(w_{n-\Delta-1},w_2) &\cdots & \alpha Tr(w_{n-\Delta-1})
\end{pmatrix}
$. It can be seen that the quotient matrix $B$ of the matrix $D_{\alpha}(G)$ is
\begin{align*} B =
\begin{pmatrix}
\frac{2\Delta^2(1-\alpha)+\alpha(t_{v_i}\Delta+Tr(v_i))}{\Delta+1}& \frac{(1-\alpha)(t_{v_i}\Delta+Tr(v_i)-2\Delta^2)}{\Delta+1}\\
\frac{(1-\alpha)(t_{v_i}\Delta+Tr(v_i)-2\Delta^2)}{n-\Delta-1} & \frac{\alpha(2W-t_{v_i}\Delta-Tr(v_i))+(1-\alpha)(2W-2t_{v_i}\Delta-2Tr(v_i)+2\Delta^2)}{n-\Delta-1}
\end{pmatrix}.
\end{align*}
The eigenvalues of the matrix $B$ are
\begin{align*}
x_1,x_2=\frac{\alpha_i\pm \sqrt{\alpha_{i}^2-4\beta_i(\Delta+1)(n-\Delta-1)}}{2(\Delta+1)(n-\Delta-1)},
\end{align*} where $\alpha_i=\alpha n(t_{v_i}\Delta+Tr(v_i)-2\Delta^2)+2n\Delta^2+(\Delta+1)(2W-2t_{v_i}\Delta-2Tr(v_i))$ and $\beta_i=2\alpha W(t_{v_i}\Delta+Tr(v_i)-2\Delta^2)+4W\Delta^2-(t_{v_i}\Delta+Tr(v_i))^2$.
Now, using Lemma 3.1, the result follows in this case.\\
\indent If $\Delta=n-1$, then $G$ being a bipartite graph implies that $G\cong K_{1,n-1}$. The result now follows from Lemma 3.3, by taking $r=1$ and $s=n-1$.\qed

\indent Taking $r=a$ and $s=n-a$ with $s\geq r$ in Lemma 3.4, it follows that the generalized distance spectrum of $K_{a,n-a}$ is
\begin{align*}
\{\alpha(n+a)-2^{[a-1]},\alpha(2n-a)-2^{[s-1]},\frac{(\alpha+2)n-4\pm \sqrt{\sigma}}{2}\},
\end{align*}
where $\sigma=n^2\alpha^2-(n^2+2a^2-2an)4\alpha+4(n^2-3an+3a^2)$. We have
$\frac{(\alpha+2)n-4+ \sqrt{\sigma}}{2}\geq \alpha(2n-a)-2$ giving  $n(2-3\alpha)+2a\alpha+\sqrt{\sigma}\geq 0$. Consider the function $f(\alpha)=n(2-3\alpha)+2a\alpha+\sqrt{\sigma},~~~~\alpha \in [0,1]$. It is easy to see that $f(\alpha)$ is decreasing function of $\alpha$. Therefore, $f(\alpha)\geq f(1)=0$ implies that $n(2-3\alpha)+2a\alpha+\sqrt{\sigma}\geq 0$ holds for all $\alpha$. Thus, it follows that $\partial_{1}(K_{a,n-a})=\frac{(\alpha+2)n-4+\sqrt{\sigma}}{2}$. Proceeding, similarly it can be seen that $\frac{(\alpha+2)n-4\pm \sqrt{\sigma}}{2}\geq \alpha(n+a)-2$ for all $\alpha\in [0,1]$ and $a\ne 1$. Therefore, for $a\ne 1$, we have $\partial_{n}(K_{a,n-a})=\alpha(n+a)-2$.\\
\indent If $\alpha\ne 0$ and $a\ne 1$, then $\partial_{n}(K_{a,n-a})=\alpha(n+a)-2$ and so
\begin{align*}
D_{\alpha}S(K_{a,n-a})=\frac{n(2-\alpha)-2a\alpha+\sqrt{\sigma}}{2}.
\end{align*}
It is easy to see that $D_{\alpha}S(K_{a,n-a})$ is a decreasing function of $a$ for all $a\in [1,\frac{n}{2}]$. Therefore, it follows that
\begin{align*}
D_{\alpha}S(K_{2,n-2})\geq D_{\alpha}S(K_{3,n-3})\geq \cdots\geq D_{\alpha}S(K_{\lfloor\frac{n}{2}\rfloor,\lceil\frac{n}{2}\rceil}).
\end{align*}
If $\alpha\ne 0$ and $a= 1$, then $\partial_{n}(K_{1,n-1})=\frac{(\alpha+2)n-4-\sqrt{\sigma}}{2}$ and so
\begin{align*}
D_{\alpha}S(K_{1,n-1})=\sqrt{n^2\alpha^2-(n^2+2-2n)4\alpha+4(n^2-3n+3)}.
\end{align*}
By direct computation it can be seen that $$\frac{n(2-\alpha)-4\alpha+\sqrt{\sigma}}{2}\leq \sqrt{n^2\alpha^2-(n^2+2-2n)4\alpha+4(n^2-3n+3)},$$ giving  $D_{\alpha}S(K_{1,n-1})\geq D_{\alpha}S(K_{2,n-2})$
Thus,  for $\alpha\ne 0$, it follows that $K_{\lfloor\frac{n}{2}\rfloor,\lceil\frac{n}{2}\rceil}$ is the complete bipartite graph of order $n$ having minimum generalized distance spectral spread and $K_{1,n-1}$ has the maximum generalized distance spectral spread.\\
\indent If $\alpha=0$, proceeding similarly as above it can be seen that the same conclusion holds. Thus we have proved the following.
\begin{theorem}
Among all the complete bipartite graphs $K_{a,n-a}$, $a\leq n-a$ of order $n$, $K_{\lfloor\frac{n}{2}\rfloor,\lceil\frac{n}{2}\rceil}$ has the minimal and $K_{1,n-1}$ has the maximal generalized distance spectral spread.
\end{theorem}

\indent The following observation follows from Lemma 2.2 and can be found in \cite{CHT}.
\begin{lemma}\label{xc}
Let $G$ be a connected graph of order $n$ and $\frac{1}{2}\leq \alpha \leq 1$. If $G^{'}=G-e$ is a connected graph obtained from $G$ by deleting an edge $e$, then  $\partial_{i}(G^{'}) \geq \partial_{i}(G)$, for all $1\leq i\leq n$.
\end{lemma}

\indent The following result gives a lower bound for $D_{\alpha}S(G)$ in terms of the order of the bipartite graph $G$.
\begin{theorem}
Let $G$ be a connected bipartite graph with $n\geq 3$ vertices. Then, for $\alpha=0$,
\begin{align*}
D_{\alpha}S(G)\geq n+\sqrt{\lceil\frac{n}{2}\rceil^2+\lfloor\frac{n}{2}\rfloor^2-\lfloor\frac{n}{2}\rfloor \lceil\frac{n}{2}\rceil},
\end{align*} with equality if and only if $G\cong K_{\lfloor\frac{n}{2}\rfloor,\lceil\frac{n}{2}\rceil}$.
For $\frac{1}{2}\leq \alpha \leq 1$,
\begin{align*}
D_{\alpha}S(G)\geq \frac{\alpha (n-3)+2n-6+ \Theta}{2},
\end{align*} where $\Theta=\sqrt{n^2\alpha^2-4(\alpha-1)(\lfloor\frac{n}{2}\rfloor^2+\lceil\frac{n}{2}\rceil^2)-4\lfloor\frac{n}{2}\rfloor \lceil\frac{n}{2}\rceil}+\sqrt{9\alpha^2-20\alpha+12}$,
equality occurs if and if only $G\cong K_{1,2}$.
\end{theorem}
{\bf Proof.} If $\alpha=0$, the proof follows from Theorem 3.4 of \cite{yzlws}. Suppose that $\alpha\ne 0$.
 Let $V(G)$ be the vertex set of the bipartite graph $G$ and let $V(G)=V_1\cup V_2$ be the bipartite partition of $V(G)$. Let $|V_1|=r$ and $|V_2|=s$, with $r+s=n$. Clearly, $G$ can be obtained by deleting some edges in $K_{r,s}$. Therefore, for $\frac{1}{2}\leq \alpha \leq 1$, from Lemma \ref{xc} and Lemma 3.4, it follows that $\partial_{1}(G) \geq \partial_{1}(K_{r,s})=\frac{\alpha n+2n-4+ \sqrt{(r^2+s^2)(\alpha-2)^2+2rs(\alpha^2-2)}}{2}$. As seen above that the quantity  $\partial_{1}(K_{r,s})=\frac{\alpha n+2n-4+ \sqrt{(r^2+s^2)(\alpha-2)^2+2rs(\alpha^2-2)}}{2}$ is minimal for $r=\lfloor\frac{n}{2}\rfloor$ and $s=\lceil\frac{n}{2}\rceil$. Thus, it follows that $\partial_{1}(G) \geq \frac{\alpha n+2n-4+ \sqrt{(\lfloor\frac{n}{2}\rfloor^2+\lceil\frac{n}{2}\rceil^2)(\alpha-2)^2+2\lfloor\frac{n}{2}\rfloor \lceil\frac{n}{2}\rceil(\alpha^2-2)}}{2}$. Since $n\geq 3$, it is clear that $K_{1,2}$ is a subgraph of $G$ and so it follows that $D_{\alpha}(K_{1,2})$ is a principle matrix of $D_{\alpha}(G)$. By Lemma 3.4, we have $\partial_{3}(K_{1,2})=\frac{3\alpha+2-\sqrt{9\alpha^2-20\alpha+12}}{2}$ and by Lemma 3.2, we have $\partial_{n}(G)\leq \partial_{3}(K_{1,2})$. The result now follows. \qed

From Theorem 3.8, for $\alpha=0$, it follows that the graph $K_{\lfloor\frac{n}{2}\rfloor,\lceil\frac{n}{2}\rceil}$ has the minimal generalized distance spectral spread. For $\alpha\ne 0$, Theorem 3.6 gives an insight that this may be true, in general. Therefore, we leave the following problem.
\begin{problem}
Among all bipartite graphs of order $n$, the graph $K_{\lfloor\frac{n}{2}\rfloor,\lceil\frac{n}{2}\rceil}$ has the minimal generalized distance spectral spread.
\end{problem}

\section{\bf Bounds for $D_{\alpha}S(G)$ in terms of clique number and independence number}

In a graph $G$, a clique is a maximal complete subgraph and the order of the maximal clique is called the clique number of the graph $G$ and is denoted by $\omega(G)$ or $\omega$. In this section, we obtain  lower bounds for generalized distance spectral spread $D_{\alpha}S(G)$ in terms of the clique number and the independence number of the graph $G$.\\

\indent First we obtain a lower bound for $D_{\alpha}S(G)$ in terms of clique number $\omega$ and the Weiner index $W$ of the graph $G$.

\begin{theorem}
Let $G$ be a connected graph on $n\geq 3$ vertices having clique number $\omega\geq 2$
 and Wiener index $W$. Let $G_1,G_2,\dots,G_k$ be all cliques of order $\omega$ and $s_i=\sum\limits_{v_j\in V(G_i)}Tr(v_j)$. \\
{\bf (i)}. If $\omega \leq n-1$, then
\begin{align*}
D_{\alpha}S(G)\geq \max_{1\leq i\leq k}\frac{\sqrt{\alpha_i^{2}-4\beta_i\omega(n-\omega)}}{(\omega)(n-\omega)},
\end{align*}  where $\alpha_i=s_i(\alpha n-2\omega)+\omega(1-\alpha)(2W+n(\omega-1))$ and $\beta_i=2 W\omega(\omega-1)-s_i^{2}+2W\alpha(s_i-\omega(\omega-1))$.\\
{\bf (ii)}. If $\omega= n$, then $D_{\alpha}S(G)=(1-\alpha)n.$
\end{theorem}
{\bf Proof.} Let $G$ be a connected  graph on $n\geq 3$ vertices having clique number  $\omega$ and  Wiener index $W$. Let $G_i$, $1\leq i\leq k$ be the cliques of $G$ having clique number $\omega$.\\
{\bf (i).} Let  $\omega\leq n-1$. For $1\leq i\leq k$, let $V(G_i)=\{u_{i_1},u_{i_2},\dots,u_{i_{\omega}}\}$ be the vertex set of $G_i$. Let $V(G)\setminus V(G_i)=\{w_1,w_2,\dots,w_{n-\omega}\}$. By suitably labelling the vertices of $G$, it can be seen that the generalised distance matrix of the graph $G$ can be written as
$D_{\alpha}(G) =
\begin{pmatrix}\label{y}
X & (1-\alpha) Y\\
(1-\alpha) Y^t & Z
\end{pmatrix},$ where
\begin{align*}
X =\begin{pmatrix}
\alpha Tr(u_{i_1}) & 1-\alpha& 1-\alpha&\cdots & 1-\alpha\\
1-\alpha & \alpha Tr(u_{i_2}) & (1-\alpha) & \cdots & (1-\alpha)\\
\vdots & \vdots & \vdots & \cdots & \vdots\\
1-\alpha & (1-\alpha)  & (1-\alpha) & \cdots &\alpha Tr(u_{i_\omega})
\end{pmatrix},
Y =\begin{pmatrix}
d(u_{i_1},w_1)  &\cdots & d(u_{i_1},w_{n-\omega})\\
d(u_{i_2},w_1) & \cdots & d(u_{i_2}w_{n-\omega})\\
\vdots  & \cdots & \vdots\\
 d(u_{i_{\omega}},w_1)&\cdots & d(u_{i_{\omega}}w_{n-\omega})
\end{pmatrix}
\end{align*}
and $Z=\begin{pmatrix}
\alpha Tr(w_1)& (1-\alpha)d(w_1,w_2)  &\cdots & (1-\alpha)d(w_1,w_{n-\omega})\\
(1-\alpha)d(w_2,w_1)& \alpha Tr(w_2)  & \cdots & (1-\alpha)d(w_2,w_{n-\omega})\\
\vdots  & \vdots& \cdots & \vdots\\
 (1-\alpha)d(w_{n-\omega},w_1)& (1-\alpha)d(w_{n-\omega},w_2) &\cdots & \alpha Tr(w_{n-\omega})
\end{pmatrix}
$.

It can be seen that the quotient matrix $B$ of the matrix $D_{\alpha}(G)$ is
\begin{align*} B =
\begin{pmatrix}
\frac{(1-\alpha)\omega(\omega-1)+\alpha s_i}{\omega}& \frac{(1-\alpha)(s_i-\omega(\omega-1))}{\omega}\\
\frac{(1-\alpha)(s_i-\omega(\omega-1))}{n-\omega}& \frac{2W+s_i(\alpha-2)+\omega(1-\alpha)(\omega-1)}{n-\omega}
\end{pmatrix}.
\end{align*}
The eigenvalues of the matrix $B$ are
\begin{align*}
x_1,x_2=\frac{\alpha_i\pm \sqrt{\alpha_{i}^2-4\beta_i\omega(n-\omega)}}{2\omega(n-\omega)},
\end{align*} where $\alpha_i=s_i(\alpha n-2\omega)+\omega(1-\alpha)(2W+n(\omega-1))$ and $\beta_i=2 W\omega(\omega-1)-s_i^{2}+2W\alpha(s_i-\omega(\omega-1))$.
Now, using Lemma 3.1, the result follows in this case.\\
\indent If $\omega=n$, then $G$ is $ K_{n}$. Since the generalized distance spectrum of the complete graph $K_n$ is $\{n-1,n\alpha-1^{[n-1]}\}$, the result follows in this case also.\qed

 A complete split graph, denoted by $CS_{t,n-t}$, is the graph consisting of a clique on $t$ vertices and an independent set on the remaining $n-t$ vertices, such that each vertex of the clique is adjacent to every vertex of the independent set. The following observation follows from Lemma 3.3 and can be found in \cite{CHT}.
 \begin{lemma}
 The generalized distance spectrum of the complete split graph $CS_{t,n-t}$ is
 \begin{align*}
 \{\alpha n-1^{[t-1]}, \alpha(2n-t)-2^{[n-t-1]}, x_1,x_2\},
 \end{align*} where $x_1,x_2=\frac{2n-t+\alpha n-3\pm \sqrt{\theta}}{2}$, $\theta=(5-4\alpha)t^2+(6\alpha n-8n-4\alpha+6)t+n^2(\alpha-2)^2+2n\alpha-4n+1$.
 \end{lemma}

\indent If $G$ is a connected graph with independence number $t=1$, then clearly $G\cong K_n$ and so $D_{\alpha}S(G)=(1-\alpha)n.$ Therefore, we need to consider $t\geq 2$. The next theorem gives a lower bound for $D_{\alpha}S(G)$ in terms of independence number $t$ of the graph $G$.
\begin{theorem}
Let $G$ be a connected  graph with $n\geq 3$ vertices having independence number $t\geq 2$. Then for $\alpha=0$,
\begin{align*}
D_{\alpha}S(G)\geq \frac{n+t+1}{2}+\sqrt{(n-t+1)^2+4t^2-4t},
\end{align*} with equality if and only if $G\cong CS_{t,n-t} $.
For $\frac{1}{2}\leq \alpha \leq 1$,
\begin{align*}
D_{\alpha}S(G)\geq \frac{2n-t+\alpha (n-3)-5+ \sqrt{\theta}+\sqrt{9\alpha^2-20\alpha+12}}{2},
\end{align*} where $\theta=(5-4\alpha)t^2+(6\alpha n-8n-4\alpha+6)t+n^2(\alpha-2)^2+2n\alpha-4n+1$,
equality occurs if and if only $G\cong K_{1,2}$.
\end{theorem}
{\bf Proof.} If $\alpha=0$, the proof follows from Theorem 3.6 of \cite{yzlws}. Suppose that $\alpha\ne 0$.
 Let $V(G)$ be the vertex set of the graph $G$ and let $V(G)=V_1\cup V_2$ be the  partition of $V(G)$ such that the subgraph induced by $V_1$ is a clique and the subgraph induced by $V_2$ is an empty graph. Let $|V_1|=n-t$ and $|V_2|=t$. Clearly, $G$ can be obtained by deleting some edges in the complete split graph $CS_{t,n-t}$. Therefore, for $\frac{1}{2}\leq \alpha \leq 1$, from Lemma \ref{xc} and Lemma 3.4, it follows that $\partial_{1}(G) \geq \partial_{1}(CS_{t,n-t})=\frac{2n-t+\alpha n-3+ \sqrt{\theta}}{2}$, where $\theta$ is defined above.  Since $n\geq 3$, it is clear that $K_{1,2}$ is a subgraph of $G$ and so it follows that $D_{\alpha}(K_{1,2})$ is a principle matrix of $D_{\alpha}(G)$. By Lemma 3.4, we have $\partial_{3}(K_{1,2})=\frac{3\alpha+2-\sqrt{9\alpha^2-20\alpha+12}}{2}$ and by Lemma 3.2, we have $\partial_{n}(G)\leq \partial_{3}(K_{1,2})$. The result now follows. \qed

{\noindent \bf Acknowledgements:} The research of S. Pirzada is supported by SERB-DST, New Delhi under the research project number MTR/2017/000084 and A. Alhevaz was in part supported by the grant from Shahrood University of Technology, Iran.

\end{document}